\documentclass[preprint,12pt]{elsarticle}
\usepackage{amssymb}
\newcounter{bla}

\usepackage{blindtext}
\usepackage[utf8]{inputenc}
\usepackage{amsmath}
\usepackage{amsfonts}
\usepackage{amssymb}
\usepackage{graphicx}
\usepackage{caption}
\usepackage{subcaption}
\usepackage{stackengine}
\usepackage{multirow}
\usepackage{array}
\usepackage{booktabs, bigdelim, rotating}
\usepackage{textcomp}
\newcommand{\textapprox}{\raisebox{0.5ex}{\texttildelow}}

\journal{Journal of Computational and Applied Mathematics}

\begin{document}

\begin{frontmatter}



\title{{~ \hfill {\small \it IMSc/2023/04}} \\ [1.5cm]
Blockwise inversion and algorithms for inverting large partitioned matrices}


\author[a,b]{R. Thiru Senthil\corref{author}}

\cortext[author] {\textit{E-mail address:} rtsenthil@imsc.res.in, ORCiD: 0000-0002-7142-4527}
\address[a]{Homi Bhabha National Institute, Anushakti Nagar, Mumbai 400094, India.}
\address[b]{The Institute of Mathematical Sciences, Taramani, Chennai 600113, India.}

\begin{abstract}
Block matrix structure is commonly arising is various physics and engineering applications.  There are various advantages in preserving the blocks structure while computing the inversion of such partitioned matrices.  In this context, using the blockwise matrix inversion technique, inversions of large matrices with different ways of memory handling are presented, in this article.  An algorithm for performing inversion of a matrix which
is partitioned into a large number of blocks is presented, in which
inversions and multiplications involving the blocks are carried out
with parallel processing. Optimized memory handling and
efficient methods for intermediate multiplications among the partitioned
blocks are implemented in this algorithm.  The developed programs for the procedures discussed in this article are provided in C language and the parallel processing methodology is implemented using OpenMP application programming interface.  The performance and the advantages of the developed algorithms are highlighted.
\end{abstract}

\end{frontmatter}

\section{Introduction}

The requirement of matrix inversion for solving a system of linear
equations arises in countless scientific and engineering applications.  Among the various methods available for computing matrix inversion, blockwise inversion methodology which computes the inversion by partitioning the matrix in blocks form, is explored in this article in programming perspective for various partitioning schemes.  For matrices of large order, block matrix inversions are performed iteratively and we discuss the possibility of efficient memory handling during the iterative processes. We also developed an
algorithm for inverting a matrix which is partitioned as a large number of
blocks. This method can be in principle applicable for inverting a general
matrix of large order, which is blocked such that the diagonal blocks are
square matrices and as an added advantage, inversions of diagonal blocks
and multiplications among the partitioned blocks in the intermediate
calculations can be carried out using parallel processing. We discuss
the development of this algorithm in this article.

General purpose programs for the procedures which are discussed in the
article including program to perform in-place inversion, are written in
C language and provided with this article. The program for the algorithm
developed in this article for inverting large partitioned matrix using
parallel processing was written using OpenMP Application Programming
Interface (API) in C language.  This project is named as {\sf `invertor'} and the written codes are available in the web-page {\sf https://www.imsc.res.in/\textapprox rtsenthil/invertor.html} and {\sf https://github.com/rthirusenthil/invertor}.

\section{Blocks partitioned matrices}

The block matrix structures are appearing naturally in various physical systems.  Some problems in fluid mechanics~\cite{fluid1, fluid2}, electrical engineering~\cite{elect} require the treatment of block partitioned matrices.  The block structure can also be seen in the systems described with multiple variables and in quantum physics~\cite{quantum1}.  The block matrix representation is often seen in density matrices of various quantum systems~\cite{quantum2}.  By preserving the block structure while performing the computations for calculating such as inversion, exponential, logarithm and determinant of the matrices, several advantages can be obtained.  One is the reduction of computational burden for calculating those matrix functions.  The other advantage is meaningful and efficient comparison of the computed matrices to the original, since the block structure is usually preserved in the final output at the end of the calculations.  Sometimes, the functions are previously known to the partitioned blocks which could be utilized to further calculate the same functions for the complete matrix structure.  For example, the work that carried out in this article address the situation that the inversion of diagonal blocks in block partitioned matrix is known and they are combined by performing multiplication and reduction while retaining the blocks structure to calculate the inversion of complete matrix.   

The block matrix structure is encountered in our studies of tau neutrino events at the Iron Calorimeter (ICAL) detector in the proposed underground India based Neutrino Observatory (INO)~\cite{taupaper}.  The pull method is used in our studies to calculate the confidence level for the appearance of tau neutrinos and their sensitivity to the neutrino oscillation parameters, with systematics.  Five types of systematics are considered in the analysis which led to the situation of inverting matrices of order $5 \times 5$ recursively to calculated minimum chi-square values, for the detector response to tau neutrinos.  The analysis required the treatment of combined samples of muon and tau type neutrinos and their anti-particles and the computation of pull variables required inverting matrices of block structure in which the complete matrix have order of $11 \times 11$ while the diagonal blocks have the order of $5 \times 5$ and $6 \times 6$.  The first diagonal block have the entries corresponding to the systematics of neutrinos and the other block have corresponding entries to the five systematics of anti-neutrinos with additional term for the ratio of number events with respect to charged flavor ratio.  It is to be noted that the non-diagonal blocks remained as non-square matrices in the block structure and have entries which are non-zero due to correlation among the various systematic variables.  The calculation of inversion of such block structured matrices lead us to develop the standalone codes for inverting block structure matrices.  The programs which are presented in the first part of this article corresponds to various possibilities of inverting the block-partitioned matrices and written in general purpose way, by expanding our blockwise inversion technique which is used in the tau analysis.  The programs address various block matrix inversion scenarios arising in different physics applications.

In data analysis, the correlation matrix that provide the relationship between the variables sometimes naturally occurs or could be arranged, in the block structure form~\cite{correlation}.  Further, the estimation of statistical significance of the model to the observation can be greatly simplified if such a decomposition is used.  Canonical representation for such block correlation matrices is discussed in~\cite{canonical}. 

The dynamical behavior of several natural systems have network descriptions and the modular random network are often encountered.  The emergence of modular structure which is the block partitioned matrix structure in the cortico-cortical networks in the brain, is shown in \cite{sitabhra}.  The modular structure behavior is common in nature for certain class of networks analysis called small-world networks (SWN) which are appearing in the studies of human society, finance, cellular metabolism and other research areas~\cite{modular}.

Thus, we have an understanding that there is requirement for treatment of such matrices which are partitioned in to large number of blocks in rows and columns, across many research areas.  In the second part of this article, we develop a general algorithm to address the scenario of calculating inverses of such large block partitioned matrices.
  
\section{Block matrix inversion}
\label{sec:Block}

A non-singular matrix ($X$) which needs to be inverted, is partitioned
into $2 \times 2$ block form as \cite{book1, book2},
\begin{equation}
X_{(m_x \times m_x)}= \left[ \begin{matrix}
  A_{(m_A \times m_A)} & B_{(m_A \times m_D)} \\
  C_{(m_D \times m_A)} & D_{(m_D \times m_D)}
\end{matrix} \right]~.
\label{eq:matx}
\end{equation}
In this form, blocks $A$ and $D$ are square matrices of order $(m_A \times m_A)$ and
$(m_D \times m_D)$ respectively, while blocks $B$ and $C$ are not necessarily square matrices.  The order of the input square matrix and the partitioned blocks are
mentioned in parenthesis at their subscripts. Here, the choice of the
order of the diagonal blocks is arbitrary. Once the orders of diagonal
blocks $A$ and $D$ are fixed, the order of $B$ and $D$ blocks are also
settled as given in the Eq.~\ref{eq:matx}.

If the block $A$ and its Schur complement $S_A=D-CA^{-1}B$ are invertible
for the partitioned matrix as in Eq.~\ref{eq:matx}, then the inversion
for matrix $X$ is given by,
\begin{equation}
X^{-1} = \left[ \begin{matrix}
A^{-1}+A^{-1}BS_A^{-1}CA^{-1} & -A^{-1}BS_A^{-1} \\
  -S_A^{-1}CA^{-1} & S_A^{-1}
\end{matrix} \right]~.
\label{eq:a}
\end{equation}

If the block $D$ and its Schur complement $S_D=A-BD^{-1}C$ are invertible, then the inverse of the matrix is given by,
\begin{equation}
X^{-1} = \left[ \begin{matrix}
S_D^{-1} & -S_D^{-1}BD^{-1} \\
  -D^{-1}CS_D^{-1} & D^{-1}+D^{-1}CS_D^{-1}BD^{-1}
\end{matrix} \right]~.
\label{eq:d}
\end{equation}

If block $A$ or its Schur's complement $S_A$ is not invertible, then the blockwise inversion can be tried using Eq.~\ref{eq:d} where blocks $D$ and $S_D$ are invertible.  In this way of blockwise inversion, the choice for order of the blocks
could be varied by keeping $m_A + m_D = m_X$ as fixed. If a
particular choice of blocking the matrix does not yield suitable non
singular blocks whose inverses are required to calculate the complete
matrix inversion, then the order of diagonal blocks among each other could be varied in the partitioning scheme, to try the blockwise inversion.

There are non singular matrices such as permutation matrices where any choice for the order of blocks $A$ and $D$ will result in situation with determinant of $A$ and $D$ or their Schur complement's determinant as zero \cite{twobytwopermute}. But those matrices might be invertible if partitioning is performed by keeping blocks $B$ and $C$ as necessarily square matrices as given by,
\begin{equation}
X_{(m_x \times n_x)}= \left[ \begin{matrix}
  A_{(m_B \times m_C)} & B_{(m_B \times m_B)} \\
  C_{(m_C \times m_C)} & D_{(m_C \times m_B)}
\end{matrix} \right]~.
\label{eq:matx1}
\end{equation}

If the block $B$ and its Schur complement $S_B = C-DB^{-1}A$ are invertible for the block matrix given in Eq.~\ref{eq:matx1}, then the inverse of the matrix is,
\begin{equation}
X^{-1} = \left[ \begin{matrix}
-S_B^{-1}DB^{-1} & S_B^{-1} \\
B^{-1}+B^{-1}AS_B^{-1}DB^{-1} & -B^{-1}AS_B^{-1}
\end{matrix} \right]~.
\label{eq:b}
\end{equation}
If the block $C$ and its Schur complement $S_C = B-AC^{-1}D$ are invertible for the block matrix given in Eq.~\ref{eq:matx1}, then the inverse of the matrix is,
\begin{equation}
X^{-1} = \left[ \begin{matrix}
-C^{-1}DS_C^{-1} & C^{-1}+C^{-1}DS_C^{-1}AC^{-1} \\
S_C^{-1} & -S_C^{-1}AC^{-1}
\end{matrix} \right]~.
\label{eq:c}
\end{equation}

The amount of computation involved and the usage of memory to compute
the matrix inversion remains the same for adopting any of the Eqs.~\ref{eq:a}, \ref{eq:d}, \ref{eq:b} and ~\ref{eq:c}. But, due to the difference in accessing the memory among the different adoptions, the difference in computation run-time would arise for a random matrix which is invertible via all the four ways. Cache handling plays a significant role in the performance of program.

For the matrix which is blocked according to Eq.~\ref{eq:matx}, if both
the blocks $A$ and $D$, and their corresponding Schur complement are
invertible, then the inverted matrix can also be given as,
\begin{equation}
X^{-1} = \left[ \begin{matrix}
S_D^{-1} & -A^{-1}B S_A^{-1}\\
  -D^{-1}CS_D^{-1} & S_A^{-1}
\end{matrix} \right]~.
\label{eq:ad}
\end{equation}
Adopting Eq.~\ref{eq:ad} for matrix inversion requires additional
computation since four blocks $A$, $D$ and their Schur complements have to be calculated comparing to two block inversions required with either Eqs.~\ref{eq:a} or \ref{eq:d}. Yet, there are advantages for adopting Eq.~\ref{eq:ad} in which the sub-blocks $A$ and $D$ or their Schur complements can be inverted simultaneously using parallel processing.  Also, some of the required multiplications can be carried out simultaneously. At first glance, the computation might appear the same even with parallel processing, but the reduction and significant advantage in computations can be obtained with multiplication among the blocks when the input matrix is treated as a configuration of large set of blocks. This procedure of inverting matrix by adopting Eq.~\ref{eq:ad} can be extended recursively to different order of matrices. A novel algorithm for handling such type of blocked matrix which are having
large number of blocks is developed and presented in this article.

If a matrix is partitioned with non-principal diagonal blocks as square matrix as given by Eq.~\ref{eq:matx1} and the blocks $B$, $C$ and their corresponding
Schur complements are invertible, then the inverted matrix is,
\begin{equation}
X^{-1} = \left[ \begin{matrix}
-S_B^{-1}DB^{-1} & S_B^{-1} \\
S_C^{-1} & -S_C^{-1}AC^{-1}
\end{matrix} \right]~.
\label{eq:bc}
\end{equation}
The amount of computation and memory usage when adopting Eq.~\ref{eq:bc}
is identical to the procedure of using Eq.~\ref{eq:ad} to compute
the inverted matrix. In this article we have developed the algorithm by
assuming that Eq.~\ref{eq:ad} is suitable for inverting a given general matrix. Our intention is to present the detailed procedure to tackle matrices of large order by
blocking them in large number of blocks. A similar identical procedure
using Eq.~\ref{eq:bc}, for inverting matrices which are not invertible
by the procedure presented here, is also possible, but we limit our discussions
here only to Eq.~\ref{eq:ad} in developing the algorithm because of the obvious reason of same amount of calculations and memory usage.

\section{Different procedures adopting block matrix inversion}
\label{sec:procs}

Large matrix inversions can be performed using blockwise inversion with
recursion \cite{ultra} technique. By adopting different procedures,
we have written C based programs to invert general matrix by
treating it as block matrix with the guidance of
Eqs.~\ref{eq:a}--\ref{eq:bc}. In order to compare these procedures for
their performance in computation and memory handling we discuss
procedures which adopt Eqs.~\ref{eq:a} and \ref{eq:ad}
in which the blocks present in the principal diagonal are square
matrices. The amount of computation, memory requirement involved with
Eq.~\ref{eq:c}, Eq.~\ref{eq:d} and Eq.~\ref{eq:bc} which contains the
blocks in counter diagonal (blocks which runs from the top right corner
to the bottom left of the matrix) as square matrices, are identical to
the procedures treating principal diagonal blocks as square matrices. Here we discuss the performance of the three written functions {\sf 'invertor\_by\_a', `invertor\_inplace\_by\_a'} and {\sf `invertor\_by\_ad'} in detail.

\subsection{`invertor\_by\_a' function}

\begin{center}
\begin{figure}[!htbp]
\centering \includegraphics[scale=0.77]{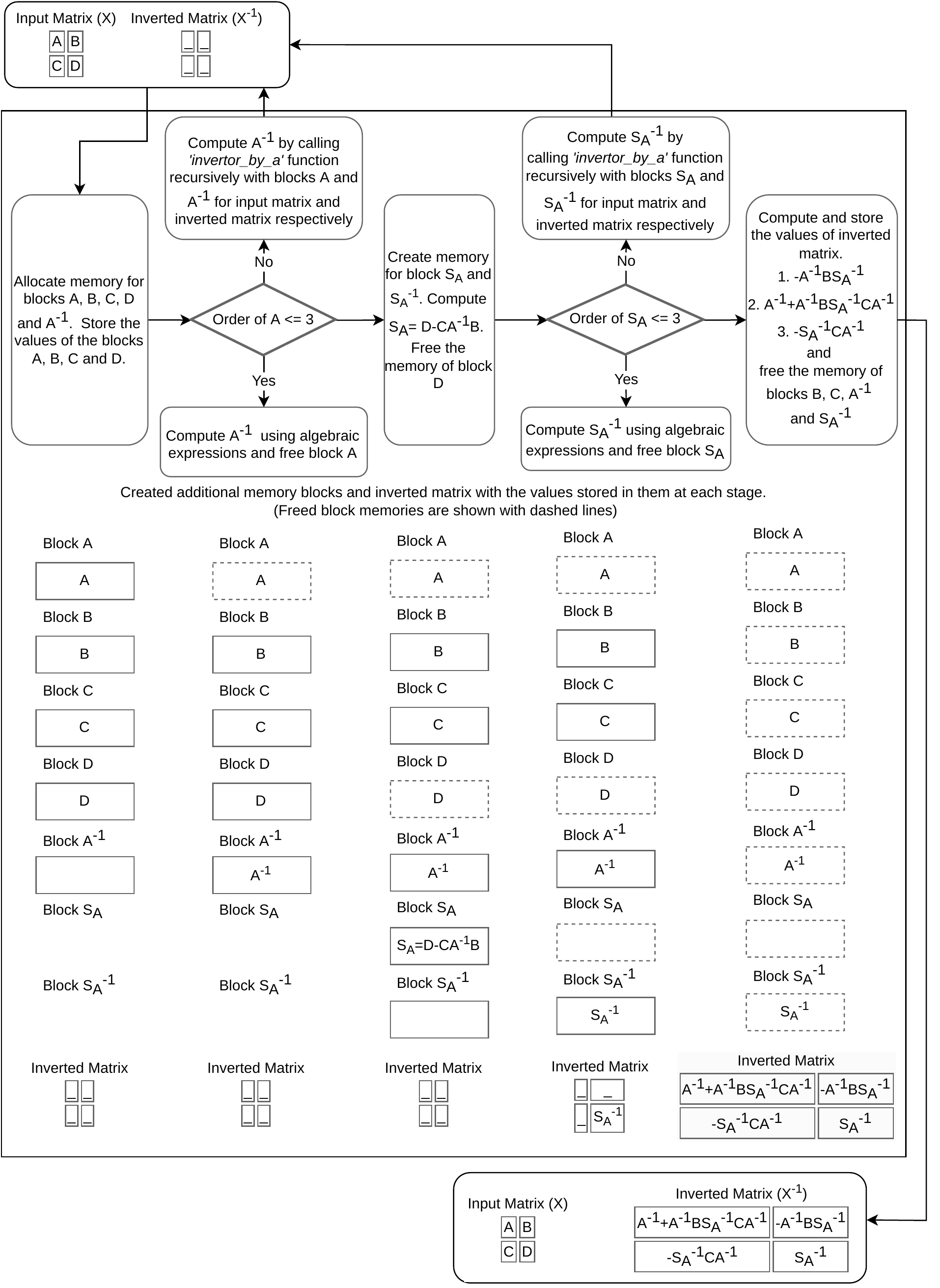}
\caption{Flowchart of {\sf `invertor\_by\_a'} function - Standard
blockwise inversion utilizing additional memory for recursive partitioned
blocks. The solid line boxes represents the active memory blocks and boxes with dashed lines indicate the memory blocks which are no longer needed and were freed, in each stage of the computation.  At the end of the procedure inverted matrix $X^{-1}$ is computed and stored at Inverted Matrix $X^{-1}$, for the given input matrix $X$.}
\label{fig:inversionfull}
\end{figure}
\end{center}

This {\sf `invertor\_by\_a'} function performs block inversion by adopting
Eq.~\ref{eq:a}.  It is the simplest approach to compute the blockwise inversion of the matrix by creating additional memory spaces for sub-blocks of blocked matrix. The computation starts with the inversion of block $A$. The flow chart of complete workflow for this process is given in Fig.~\ref{fig:inversionfull}.  
Apart from the input matrix and inverted matrix memory, the computation requires creation of block memories dynamically during the execution.  In Fig.~\ref{fig:inversionfull}, additional memory spaces used in each step are also
shown.  The block $A$ and the other created block, the Schur's complement $S_A$, must be non-singular square matrices. In this procedure, it is required to calculate the inverse of both $A$ and $S_A$: 2 inversions. The multiplications among the blocks which need to be computed are,
\begin{itemize}
\item $-C*A^{-1}*B$ - multiplications of three matrices.
\item $-A^{-1}*B*S_A^{-1}$ - multiplication of three matrices.
\item $-S_A^{-1}*C*A^{-1}$ - multiplication of three matrices.
\item $A^{-1}*B*S_A^{-1}*C*A^{-1}$ - multiplication of five matrices.
\end{itemize} 
Multiplying three matrices can be performed by multiplying two matrices at
first and then multiplying the result to the third matrix.  This way of multiplying requires additional storage as temporary block for the result of initial two matrices multiplication.  The multiplication of three matrix can also be performed by calculating the only the row or column of first two matrices multiplication each time in order to multiply with the column or row of the third matrix, to get every corresponding element in the three matrix multiplication.  Although this process requires lesser temporary memory (by re-using the same memory for subsequent calculations), it would not be an efficient way since same stored blocks are called in to cache for computation multiple times.  For comparing the different procedures, that discussed in this article, we do not include this temporary memory requirement and management for the multiplication among the blocks since they can be handled in multiple ways.  But, with parallel processing algorithm, we even take care of this memory management effectively which is discussed in the section \ref{sec:parallel}.

In the listed multiplications required for {\sf `invertor\_by\_a'} function, the multiplication of three matrices can be treated as two times of two matrices multiplications.  Similarly, multiplication of five matrices is equivalent to four times of two matrices multiplications.  In total, the set of multiplications required in this process is equivalently, nine multiplications of two matrices at a time. 

But if the process is carried out starting with $-A^{-1}B$ and $CA^{-1}$ - two multiplications of two
matrices, then its result could be used to compute the listed required
items with additional four multiplications appropriately with $B$ and
$S_A^{-1}$. In total only six multiplications of two matrices are
required in this process. Apart from the multiplications two reduction
sums are also required during the computation of Schur complement
$S_A$ and inverted block which corresponds to block $A$. The number
of additional memory blocks which are required in this computation are
also shown in the Fig.~\ref{fig:inversionfull}. Ignoring the order of
the block matrix, it can be seen that seven additional memory blocks
are required in this approach but only five memory blocks are actively
used in any step. The remaining blocks are created subsequently after
freeing the used memory blocks.

When the input matrix is very large, the inversion of sub blocks
is performed by calling the same function, recursively. In order to
have an understanding on the requirement of memory and computations,
let us calculate these requirements for a matrix of order $2^n$
where $n=1, 2, 3, ..., k,... $. For $n=1$, the inversion of the
matrix is performed in one iteration. It requires 2-inversions, 6-two
matrix multiplications, 2-reductions. The additional memory of seven
blocks are of the order of 1. For an input matrix of order $2^k$, the
inversion is calculated recursively. In total $2^k$ inversions will be
performed. The number of multiplications required for every call of the
inversion function (recursively) will be six. In total for the matrix of
order $2^k$, it can be seen that the total number of multiplications are
$6^{1+2^1+2^2+...+2^{k-1}}=6^{2^k-1}$. The multiplications involved in the
sub-calls of inversion function recursively will be among the matrices
of order half of that of the order of parent blocks. This estimation
provided here did not account for the computation due to the size of
matrix order, but only the count on the number of multiplications required
to perform the inversion. The scalar multiplication of $-1$ with the
matrices can be included while performing matrix multiplications. Hence,
we did not count this as separate multiplications. Similarly, the number
of reduction operations on account of all the sub-calls of the inversion
function will be $2^{2^k-1}$. When it comes to the memory required in
every sub call, the size of the memory of each block will be half that of
the parent function. If we assume the recursive calls for the inversion
of blocks are called till sub-block matrix of order $2\times2$, then
the total number of additional memory required for inverting matrix
of order $2^k$ will be $7*[2^{k-1}(2^0*2^0)+2^{k-2}(2^1*2^1)+ \ldots +
2^1(2^{k-2}*2^{k-2})+2^0(2^{k-1}*2^{k-1})]$. This series sum results
in $7*2^k*2^{k-1}$. The quantification that we have performed for this
inversion method on memory and amount of computation are useful to
compare this procedure with others for optimization.

\subsection{`invertor\_inplace\_by\_a' function}

\begin{figure}[tbp]
\centering \includegraphics[scale=0.85]{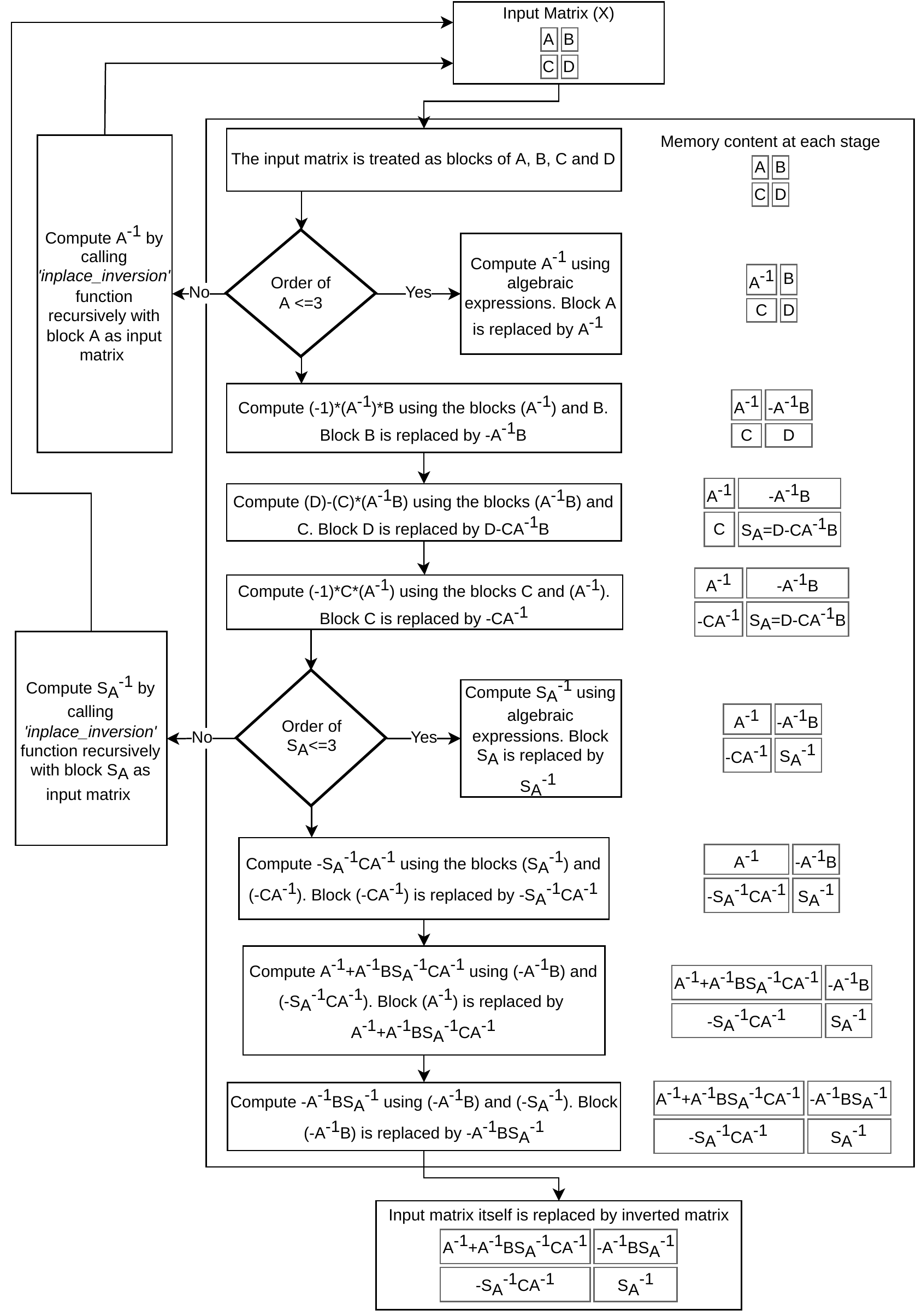}
\caption{Flowchart of {\sf `invertor\_inplace\_by\_a'} function -
Complete blockwise inversion which do not have the requirement for
additional memory for recursive partitioned blocks. The updated values in the input matrix from the calculation for inverted matrix at every step, is also shown
in the figure.  At the end of the procedure, the input matrix itself is
replaced by the inverted matrix.}
\label{fig:inplaceinversion}
\end{figure}

The flowchart of {\sf `invertor\_inplace\_by\_a'} function is provided
in Fig.~\ref{fig:inplaceinversion}.

In this method, the input matrix itself will be replaced by the inverted
matrix at the end of the procedure. This procedure has the advantage
of not using additional memory for storing the blocks in intermediate
computations while calling the inversion recursively. Apart from the
two inversions per iteration, the required multiplications among the blocks are,

\begin{itemize}
\item $-A^{-1}*B$,
\item $C*(-A^{-1}B)$,
\item $-C*A^{-1}$,
\item $S_A^{-1}*(-CA^{-1})$,
\item $(-A^{-1}B)*(-S_A^{-1}CA^{-1})$,
\item $(-A^{-1}B)*(-S_A^{-1})$.
\end{itemize}
It is to be noted that for all the listed multiplications, when performed
in the above given order, only two block matrices are multiplied at a time. In total, six two matrices multiplications are required
for the calculation of the inverse. Similar to the analysis performed
for the inverse function, we can also estimate the amount of computation
for this function for the case of input matrix of order $2^k$. When the
subsequent block inversion are computed recursively, the total number
of block matrix multiplications required would be $6^{2^k-1}$, using
the calculation given in inverse function. The number of inversions
and reductions which need to be performed for a matrix of order $2^k$
is identical to that for the inversion function which is $2^k$.

The multiplications performed in this function are carried out
in-place, which means that the block in which the multiplication is
performed (by left or right matrix multiplication with another block),
is replaced by the result. For example, block $B$ of the input matrix
is replaced by $-A^{-1}*B$ in the first multiplication. In-place matrix
multiplications cannot avoid the usage of temporary memory at least
the size of the row/column of the actual input matrix in which the
multiplication is performed. This additional memory requirement is
required in the subsequent recursive {\sf `invertor\_inplace\_by\_a'}
function calls. While the program computes the inversion serially,
temporary memory of an array of size $k$ is sufficient for input matrix of
order $(k \times k)$. The recursive calls could use the same memory for
intermediate matrix multiplications since their order would be lesser
and lesser in the subsequent inverse calls.

\subsection{`invertor\_by\_ad' function}

\begin{figure} \centering
\includegraphics[width=0.99\textwidth,height=1.2\textwidth]{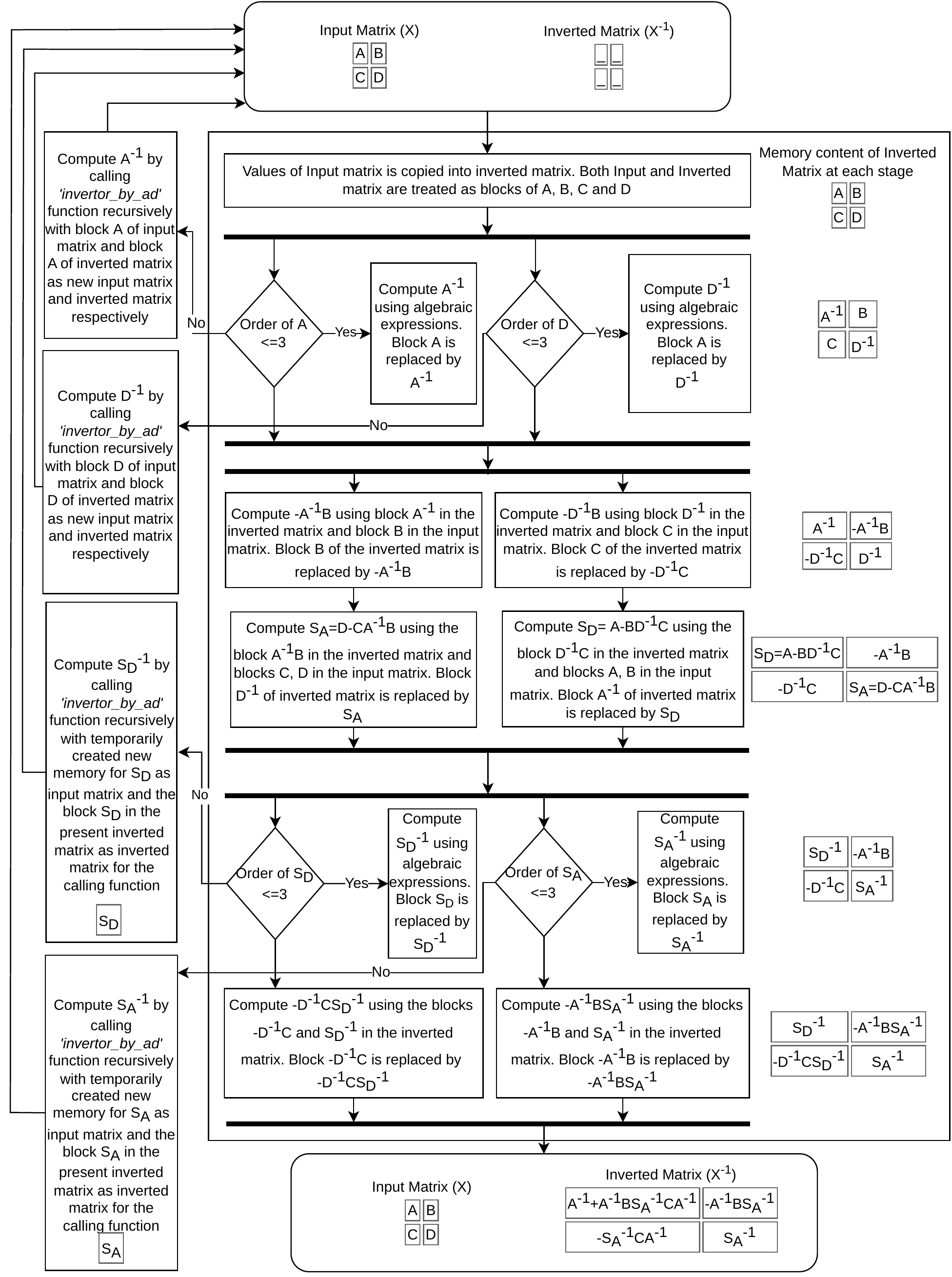}
\caption{Flowchart of {\sf
`invertor\_by\_ad'} function - The inversion and multiplications among
the blocks are carried out, in principle, simultaneously in each step.  The implicit
barriers, which are required between the steps are shown in thick solid
lines which ensure that the calculations are performed completely among the threads before proceeding to next steps if parallel processing was carried out in those steps.  Additional memory declaration required for Schur complements ($S_A$ and $S_D$) in recursive calls of the functions is also shown. The updated inverted matrix at the end of each step is also shown. At the end of the procedure inverted matrix ($X^{-1}$) is computed for the given input matrix $X$.}
\label{fig:invertbyaandd}
\end{figure}

The flowchart for this function, adopting Eq.~\ref{eq:ad} is given in
Fig.~\ref{fig:invertbyaandd}.

Here it is assumed that the blocks $A$, $D$ and their Schur complements are
invertible in the parent input matrix but also in the subsequent recursive
calls child input blocks. The inverted matrix status in each step is
also given in Fig.~\ref{fig:invertbyaandd}. The number of inverses which
need to be performed in every iteration is four, which is twice that of
inversion required with {\sf `invertor\_inplace\_by\_a'} function. The
sum reductions in this calculation are two which is the same that of {\sf
`invertor\_inplace\_by\_a'} function. We can also have an advantage of
carrying out each set of two inversions simultaneously. The order
in which the inversion is performed in each pair of inversions will not
change the result. This serves as a `guidance' for our parallel algorithm
which is discussed in detail in the next section. The disadvantage that we
see in this function is due to the requirement of additional memory for the
Schur complements which is mandatory for recursive calls. This additional
memory requirement is significantly less than what we have calculated
with {\sf `invertor\_by\_a'} function. The list of multiplications which
are needed in this {\sf 'invertor\_by\_ad'} function are:

\begin{itemize}
\item $-A^{-1}*B$,
\item $-D^{-1}*C$,
\item $(-D^{-1}C)*S_D^{-1}$,
\item $(-A^{-1}B)*S_A^{-1}$.
\end{itemize}
Hence the number of multiplications required are only four two
matrix multiplications for every call of the function, when they are
performed in the above given order. The total number of block matrix
multiplications required would be $4^{2^k-1}$ for the matrix of order
$2^k$.  Recursive calls in the function require additional memory for
the Schur complements $S_A$ and $S_D$. Additional memory for the matrix
of order $2^k$ computed for recursion taken up to the level of $2
\times 2 $ sub-blocks is $2^1(2^{k-1}*2^{k-1})+2^2(2^{k-2}*2^{k-2})+
\ldots +2^{k-1}(2^{1}*2^{1})$. The series sum results as
$2^{k+1}(2^{k-1}-1)$. This additional memory required is significantly
less than {\sf `invertor\_by\_a'} function but more than {\sf
`invertor\_inplace\_by\_a'} function. Multiplication among the blocks
would also require additional memory temporarily as estimated in {\sf
`invertor\_inplace\_by\_a'} function.

\subsection{Performance of different procedures for block inversion}

\begin{center}
\begin{figure}[htb]
\centering
\includegraphics[scale=0.48, angle=270]{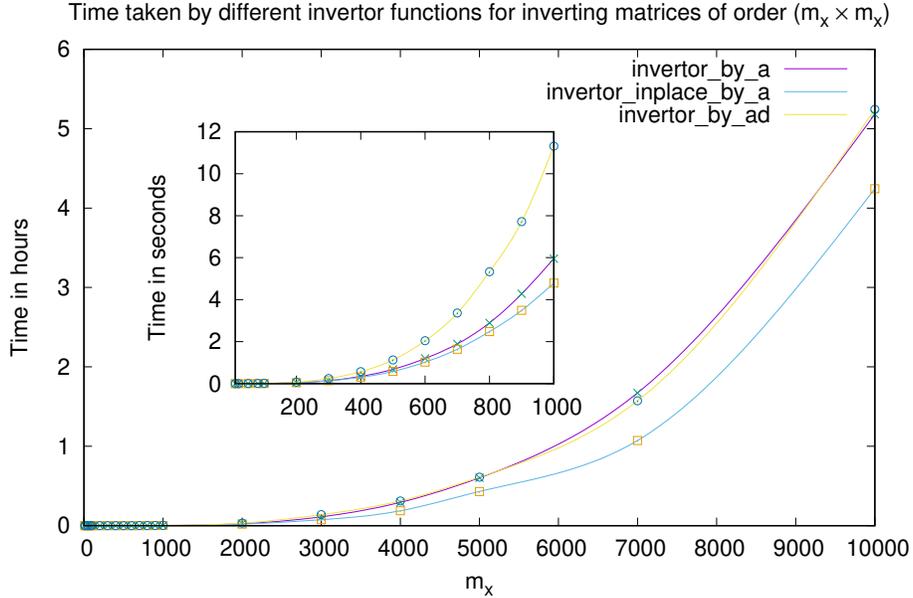}
\caption{Time taken for inverting matrices up to $m_x=10^4$ using different functions discussed in the article.  The inside figure shows the computation time for matrices of order up to 1000.}
\label{fig:invertorstime}
\end{figure}
\end{center}

The different functions developed here have different advantages.  C based programs are written and tested with a common laptop machine whose specifications are: Intel Core i5-10210U CPU 1.60GHz processor, 8GB RAM machine.  The time taken for computing inversion among the different functions are plotted in Fig.~\ref{fig:invertorstime}.  In-place inversion procedure is efficient in computation time for serial computation of matrix inversion. Yet this procedure still requires additional memory during the multiplication between the blocks which is unfriendly for cache handling when the size of the block is large. The
{\sf `invertor\_by\_a'} function requires more additional memory usage
which also results in more computation time due to data handling. But
the multiplication among the blocks can be performed in cache friendly
manner wherever possible. The {\sf `invertor\_by\_ad'} function has the
advantage of computing the inversions and multiplications simultaneously
among different stages.  Although the {\sf `invertor\_inplace\_by\_a'}
function is highly effective among the three procedures that is discussed
here, it is not readily adoptable for parallel processing, where recursion
should be avoided and the additional memory required during intermediate
computations has to be identified and declared before the execution of
parallel computations.

\begin{center}
\begin{figure}[htb]
\centering
\includegraphics[scale=0.48, angle=270]{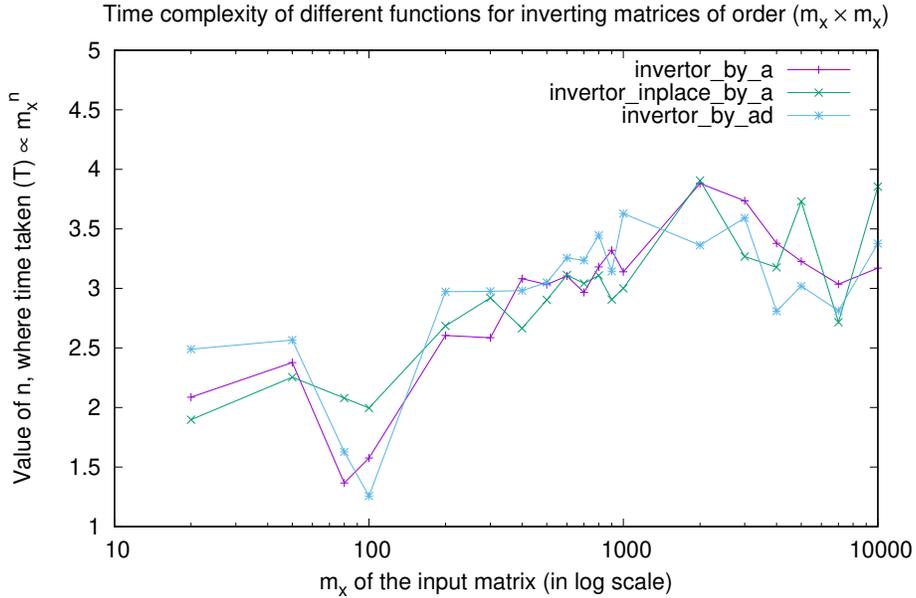}
\caption{Variation of time complexity for different regions of order of the input matrix.  Here, the time taken for computing matrix inversion ($T$) is assumed as $T \propto m_x^n$.  Value of n is calculated from the derivative of $log~T$ with respect to $log~m_x$ at different values of $m_x$.  The figure shows the time complexity $\mathcal{O}\left(m_x^n\right)$ with the value of $n$ as well below $3$, for the input matrices with $m_x$ up to $10^2$.  The increase in $n$, can be realized as the effect of cache handling with multiplications and inversions of sub blocks in very large order matrices.}
\label{fig:slope}
\end{figure}
\end{center}

In order to estimate the time complexity for comparing the performance of different functions, we fist assume the amount of time spent ($T$) performing the inversion of matrix of order $(m_x \times m_x)$ as,
$$T \propto m_x^n .$$
The time complexity $\mathcal{O}\left(m_x^n\right)$ can be realized from the value of $n$.  For different regions in values of $m_x$, the calculated value of $n$ is plotted in Fig.~\ref{fig:slope}.  For range of values in $m_x$, the calculated values of $n$ is provided in Table~\ref{tab:timeslope}.  It can be seen that the $n$ is well below $3$ for all three functions that were developed for the matrices with $m_x$ up to the order of $10^2$.  The increase in $n$ for large matrices with $m_x > 10^3$ can be realized as an effect of cache handling by the machine during the multiplications among the blocks as large sized arrays.  Although there are ways to improving the cache handling \cite{cachemul}, such as segmenting the loops and performing the calculations in efficient cache friendly manner which will improve the efficiency of the written function significantly, we have not implemented the techniques in these procedures while studying their performance.  And also these functions are developed and written to achieve the possibility of developing an algorithm for carrying out inversion of highly block partitioned large matrices, which is discussed in detail in the next section. 
\begin{center}
\begin{table}
\begin{tabular}{|c|c|c|c|}
\hline
Range of $m_x$ of the input matrix & invertor\_by\_a & invertor\_inplace\_by\_a & invertor\_by\_ad \\
\hline
$10 - 100$ & 2.45 & 2.60 & 2.77 \\
$10 - 1000$ & 3.74 & 3.31 & 3.58 \\
$10 - 10000$ & 3.15 & 3.66 & 3.25 \\
\hline
\end{tabular}
\caption{Table containing the values of $n$ for time complexity $\mathcal{O}\left(m_x^n\right)$ in different range $m_x$. The increase in $n$ is realized as the effect of cache handling during the multiplications among the partitioned blocks when the array sizes become very large.}
\label{tab:timeslope}
\end{table}
\end{center}

\section{Algorithm for inverting large partitioned matrix with parallel processing} \label{sec:parallel}
For a matrix which is partitioned into large number of blocks, the development of algorithm for performing inversion is discussed in this section. From the developed three procedures we have the guidance for handling temporary memory required for the intermediate calculations.  Also, we consider the possibility of performing sub blocks inversions simultaneously and also independent of other calculations in the procedure.  Some of the considerations for the development of the parallel processing algorithm which need to be considered throughout the calculation can be summarized as follows.

\begin{itemize}
\item The inversion of diagonal blocks and multiplications among the blocks can be performed using parallel processing.  It meant that the calculations are such that the operations on the blocks are completely need to be independent of each other.

\item The developed algorithm should avoid recursion and must precisely
predefined the required iterations beforehand for matrices of different orders. The bottom to top approach (i.e., starting with the inversion of blocks of smallest
size to the larger block till the size of the complete input matrix is
reached) is suitable.

\item The required temporary memories required for intermediate
multiplications, reductions and for storing the calculated blocks which
are required to be used in further steps, have to be precisely determined and allocated beforehand before starting the parallel processing.  Otherwise, creation and deletion of memories in thread-specific manner might result in segmentation fault with poor memory handling. The performance of complete algorithm is therefore highly depend on precise handling of additional memory requirements.

\item In general, computation and memory management using parallel
processing do not change much from serial computation, since computing
large matrix inversion highly depends on how the memory and data are
being managed rather than the complexity of computation. Hence, the
procedure should improve the efficiency in data transfer with proper
cache handling and minimum usage of additional memory for intermediate
calculations, in order to have a better performance.

\item In addition to the above, we also consider scenario of storing and retrieving values of partitioned blocks in separate files and performing the calculations in serious of steps which do not have to be performed in one go.  It meant that the algorithm attempts to create the possibility of performing inversion to any extremely large order, computations involving multiple machines without anytime frame for completing the inversion.

\end{itemize}
By considering the restrictions listed above the given algorithm here
is optimized in three different limits in the computation which are: (i)
minimum additional memory usage, (ii) efficient data handling for cache
for blockwise multiplications, and (iii) bottom to top approach algorithm
with required control on parallelism with a suitable key variable for barriers.  From the three functions created for performing inversion of matrix whose primary diagonal blocks and their Schur complements are invertible, we infer different advantageous as follows in this algorithm.
\begin{itemize}
\item From {\sf `invertor\_by\_a'} function: Precise determination of additional memory requirement are inferred from this function for the algorithm.  The temporary memory are handled by creating provisional memory blocks set given by $T$ as matrices of different suitable orders. Over different steps in the complete procedure, these provisional set of memory blocks are reused appropriately.
\item From {\sf `invertor\_inplace\_by\_a'} function:  The inverted matrix is being updated in iterative manner with inverses of blocks in this function.  We use the same way to update the inverted matrix in the parallel processing algorithm, in which, the partitioned diagonal blocks given by $\left(X^l\right)^{-1}$ in the inverted matrix are updated with inverses of lower to higher orders sequentially such as $l=1, 2, ...$.  The steps are created to perform the inversion from the inverted diagonal blocks of one order less, which also provide the choice of a situation in which inversion of the complete matrix is not performed in one go, but with breaks in sequence of steps identified by a variable called `stepid'.
\item From {\sf `invertor\_by\_ad'} function:  In this function the inversions of diagonal blocks $A$ and $D$ can be carried out simultaneously and also the inversions of their Schur complements.  The calculations among the partitioned blocks are also completely independent of each other among each stage, which meant that they can also be performed in parallel.  We adopt this procedure in the algorithm, for carrying out inverses and multiplications in and among the blocks in parallel.  The requirement of completion of all the calculations before moving to the next step is controlled by a variable called `stepid' in the algorithm.  
\end{itemize}

In the next subsection we discuss details of partitioning the input matrix
into appropriate blocks. In the subsequent subsections, we discuss the
memory allocation for intermediate calculations, notations for appropriate
computations and the parallel computing procedure itself to compute the
matrix inversion.

\subsection{Partitioning of the input matrix into blocks}
The algorithm is designed to invert the matrix with the number of
diagonal square blocks as powers of 2. Hence, the input matrix is
partitioned as $[2^n \times 2^n]$ blocks. Matrices of different order
are accommodated in the $[2^n \times 2^n]$ blocks structure with
appropriate choice for order of sub-block matrices. In order to
distinguish order of the matrix (given by number of elements in rows
$\times$ number of elements in columns) from its block partitioned
structure (given by {\sf `blockorder'} which is number of blocks in
rows $\times$ number of blocks in columns), we follow the notation of
representing the order within parenthesis and the `{\sf blockorder}'
within the square brackets for the matrices and sub-blocks. Once the
order of diagonal blocks are fixed, the orders of non diagonal blocks
are also settled as given below. The orders for the diagonal blocks are
arbitrarily chosen. The restriction is that the partitioned block matrix
must be invertible and the sum of orders of diagonal blocks has to be
equal to the order of input matrix.
\begin{equation}
M_{(m_{n} \times m_{n})} = \left[ \begin{matrix}
P^{\langle 0\rangle \langle 0\rangle }_{(m_{1} \times m_{1})} &P^{\langle 0\rangle \langle 1\rangle }_{(m_{1} \times m_{2})} ~~& \dots ~~&P^{\langle 0\rangle \langle k-1\rangle }_{(m_{1} \times m_{k})} ~\\
P^{\langle 1\rangle \langle 0\rangle }_{(m_{2} \times m_{1})} ~~&P^{\langle 1\rangle \langle 1\rangle }_{(m_{2} \times m_{2})} ~~& \dots ~~&P^{\langle 1\rangle \langle k-1\rangle }_{(m_{2} \times m_{k})} ~\\
\vdots & \vdots & \ddots ~~& \vdots \\
P^{\langle k-1\rangle \langle 0\rangle }_{(m_{k} \times m_{1})} ~~&P^{\langle k-1\rangle \langle 1\rangle }_{(m_{k} \times m_{2})} ~~& \dots ~~&P^{\langle k-1\rangle \langle k-1\rangle }_{(m_{k} \times m_{k})} ~\\
\end{matrix} \right]~; ~m_n=\sum_{i=1}^k m_i~.
\label{eq:blocking}
\end{equation}
In Eq.~\ref{eq:blocking}, the input matrix is partitioned as blocks of
order $[k \times k]$ and the indices ($0, 1, \dots, k$) that identify the
blocks are represented within angle brackets. For the program developed
in this {\sf `invertor'} project, we find that it is sufficient to choose
the orders of the individual diagonal blocks to be among 2, 3 and 4. i.e.,
$m_i \in \{2, 3, 4\}$. Then, the total number of diagonal blocks $N_k$
can be calculated as,
\begin{equation} \nonumber
N_k = 2^{[(\rm int)log_2(m_n)]-1}~,
\end{equation}
where $(m_n \times m_n)$ is the order of the input matrix and it is
also to be noted that we convert the {\it double} value of $\log_2
(m_n)$ into {\it integer}, which ensures the correct number of blocks
corresponding to our choice for orders of partitioned blocks. In
Table~\ref{tab:blocking}, blocking scheme of diagonal blocks is shown
for various orders of input matrix.
\begin{table}[htbp]
\begin{center}
\begin{tabular}{|c |c c c c c c c c c |c|} \hline
Order $m_n$& \multicolumn{9}{|c|}{Diagonal blocks and their order} &
 Total number of blocks \\ \hline 
 & & & & & & & & & & \\
2& \fbox{2} & & & & & & & & & 1\\
3& \fbox{3} & & & & & & & & & 1\\
4& \fbox{2}& \fbox{2} & & & & & & & & 2\\
5& \fbox{2}& \fbox{3} & & & & & & & & 2\\
6& \fbox{3}& \fbox{3} & & & & & & & & 2\\
7& \fbox{3}& \fbox{4} & & & & & & & & 2\\
8& \fbox{2}& \fbox{2}& \fbox{2}& \fbox{2} & & & & & & 4 \\
9& \fbox{2}& \fbox{2}& \fbox{2}& \fbox{3} & & & & & & 4 \\
\vdots & & & \vdots & & & & & & & \vdots \\
21& \fbox{2}& \fbox{2}& \fbox{2}& \fbox{3} & \fbox{3} & \fbox{3} & \fbox{3} & \fbox{3} & & 8 \\
\vdots & & & \vdots & & & & & & & \vdots \\
$m_n$ & \fbox{$m_0$} & \fbox{$m_1$} & \fbox{$m_2$}& \fbox{$m_3$} & \ldots
& \fbox{$m_i$} & \ldots&\ldots & \fbox{$m_{{N_k}}$} & $N_k=2^{[({\rm
int})\log_2 (m_n)]-1}$  \\
 & & & & & & & & & & \\ \hline
\end{tabular}
\caption{Fixing order of diagonal blocks for input matrix of different
orders while partitioning. Here the $m_i$, order of the diagonal blocks,
are fixed by the scheme given by Fig.~\ref{fig:blocking} when
{\sf blockorder} is chosen from the set $\{2, 3, 4\}$.}
\label{tab:blocking}
\end{center}
\end{table}
\begin{center}
\begin{figure}
\centering
\includegraphics[width=0.45\textwidth, height=0.55\textwidth]{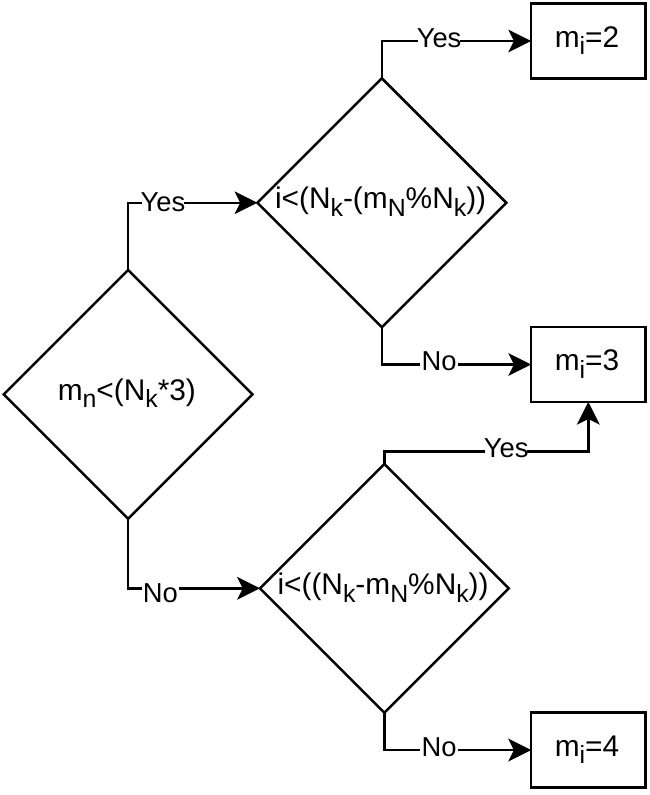}
\caption{Selecting order $m_i$ of the diagonal blocks. Here $m_i \in \{2,
3, 4\}$, $i=0, 1, 2, ... ,N_k$ and $N_k$ is the total number of diagonal
blocks given by, $N_k=2^{[({\rm int})\log_2 (m_n)]-1}$ while the order of the unpartitioned input matrix is $(m_n \times m_n)$.}
\label{fig:blocking}
\end{figure}
\end{center}
The procedure requires finding the inverses of diagonal blocks and the
Schur complements which involve blocks of larger size containing these
diagonal blocks. For instance, input matrices of order $100$ and $10000$
require $2^5$ and $2^{12}$ blocks respectively. The order of the blocks
are fixed as given in the flow chart of Fig.~\ref{fig:blocking}. Since
our choice for the order of the blocks is in the set $\{2, 3, 4\}$, the
individual diagonal block inverses can simply be computed using algebraic
expressions. If the inverse cannot be computed for any of the block matrix
which are singular with the adopted blocking scheme, then the blocks
can be rearranged with different orders by satisfying the appropriate
conditions. Even if the different choices for partitioning do not result
in non singular block matrices in all the stages of computation, the order
of blocks can be increased or varied with the restriction that the total
number of diagonal blocks must be in the powers of 2 for this algorithm.

\subsection{Calculations of inverse of input matrix diagonal blocks $X^l$}
\label{ssec:inversexl}

The matrix is partitioned into $2^k \times 2^k$ blocks, as described
in the previous section. In the bottom to top approach to invert, we
start out with inversions of the diagonal blocks. Using Eq.~\ref{eq:ad},
the inversion of blocks in which each block contains $2 \times 2$
sub-blocks for a total of $k/(2^1)$ inverted blocks, is done. This
procedure is further repeated for finding inverse of higher order
blocks in which each block contains sub-blocks of order $2^i \times 2^i,
i\in \{1,2,3,...,(k)\}$, till the inverse of the complete input matrix
is computed.

From the inverses of diagonal blocks containing $2^j \times
2^j$ blocks, to find the inverses of diagonal blocks with $2^{j+1}
\times 2^{j+1}$ blocks, the number of {\sf invert\_by\_ad} function
identical computations to be performed are $2^{(k/j)}$. The inverse for
each block is calculated by,
\begin{equation}
M_{} = \left[ \begin{matrix}
X^{\langle 1\rangle }_{[2^{j+1} \times 2^{j+1}]} & & & ~\\
& X^{\langle 2\rangle }_{[2^{j+1} \times 2^{j+1}]} & & ~\\
& & \ddots & ~\\
& & & X^{\langle 2^{k/j}\rangle }_{[2^{j+1} \times 2^{j+1}]}
\end{matrix} \right]~;
\label{eq:blocking2} 
\end{equation}
with
\begin{equation}
X^{l}_{[2^{j+1} \times 2^{j+1}]} = \left[ \begin{matrix}
A^l_{[2^{j} \times 2^{j}]} & B^l_{[2^{j} \times 2^{j}]}\\
  C^l_{[2^{j} \times 2^{j}]} & D^l_{[2^{j} \times 2^{j}]}
\end{matrix} \right],~~ l=1, 2, \dots, (2^{k/j})~.
\label{eq:adprl}
\end{equation}
Here, the order of the blocks $A^l$ is given as $m_l \times m_l$ and $l=1,2,
\ldots ,(k/j)$. For different values of $l$, the value of $m_l$ can
be computed as $m_l=\left( \sum_{i=0}^j m_{l+i}\right)$ for $A_{(m_l
\times m_l)}$. The description of blocks $A^l, B^l, C^l$ and $D^l$
in terms of partitioned blocks of the input matrix is given by,
\begin{eqnarray} \nonumber
A^l_{[2^{j} \times 2^{j}]}=&\left[ \begin{matrix}
P^{\langle l*2^j+0\rangle \langle l*2^j+0\rangle } & P^{\langle l*2^j+0\rangle \langle l*2^j+1\rangle } & \cdots & P^{\langle l*2^j+0\rangle \langle l*2^j+(2^j-1)\rangle }\\
P^{\langle l*2^j+1\rangle \langle l*2^j+0\rangle } & P^{\langle l*2^j+1\rangle \langle l*2^j+1\rangle } & \cdots &P^{\langle l*2^j+1\rangle \langle l*2^j+(2^j-1)\rangle }\\
\vdots & \vdots & \ddots &\vdots \\
P^{\langle l*2^j+(2^j-1)\rangle \langle l*2^j+0\rangle } & P^{\langle l*2^j+(2^j-1)\rangle \langle l*2^j+1\rangle } & \cdots &P^{\langle l*2^j+(2^j-1)\rangle \langle l*2^j+(2^j-1)\rangle } \\
\end{matrix} \right]~;\nonumber \\
B^l_{[2^{j} \times 2^{j}]}=&\left[ \begin{matrix}
P^{\langle l*2^j+0\rangle \langle (l+1)*2^{j}+0\rangle } & P^{\langle l*2^j+0\rangle \langle (l+1)*2^{j}+1\rangle } & \cdots &P^{\langle l*2^j+0\rangle \langle (l+1)*2^{j}+(2^j-1)\rangle }\\
P^{\langle l*2^j+1\rangle \langle (l+1)*2^{j}+0\rangle }& P^{\langle l*2^j+1\rangle \langle (l+1)*2^{j}+1\rangle }& \cdots &P^{\langle l*2^j+1\rangle \langle (l+1)*2^{j}+(2^j-1)\rangle }\\
\vdots & \vdots & \ddots &\vdots \\
P^{\langle l*2^j+(2^j-1)\rangle \langle (l+1)*2^{j}+0\rangle }& P^{\langle l*2^j+(2^j-1)\rangle \langle (l+1)*2^{j}+1\rangle }&
\cdots & P^{\langle l*2^j+(2^j-1)\rangle \langle (l+1)*2^{j}+(2^j-1)\rangle }\\
\end{matrix} \right]~; \nonumber \\
C^l_{[2^{j} \times 2^{j}]}=&\left[ \begin{matrix}
P^{\langle (l+1)*2^{j}+0\rangle \langle l*2^j+0\rangle }& P^{\langle (l+1)*2^{j}+0\rangle \langle l*2^j+1\rangle }& \cdots &P^{\langle (l+1)*2^{j}+0\rangle \langle l*2^j+(2^j-1)\rangle }\\
P^{\langle (l+1)*2^{j}+1\rangle \langle l*2^j+0\rangle }& P^{\langle (l+1)*2^{j}+1\rangle \langle l*2^j+1\rangle }& \cdots &P^{\langle (l+1)*2^{j}+1\rangle \langle l*2^j+(2^j-1)\rangle }\\
\vdots & \vdots & \ddots &\vdots \\
P^{\langle (l+1)*2^{j}+(2^j-1)\rangle \langle l*2^j+0\rangle }& P^{\langle (l+1)*2^{j}+(2^j-1)\rangle \langle l*2^j+1\rangle }&
\cdots & P^{\langle (l+1)*2^{j}+(2^j-1)\rangle \langle l*2^j+(2^j-1)\rangle }\\
\end{matrix} \right]~; \nonumber \\
D^l_{[2^{j} \times 2^{j}]}=&\left[\begin{matrix}
P^{\langle (l+1)2^{j}+0\rangle \langle (l+1)2^{j}+0\rangle }& P^{\langle (l+1)2^{j}+0\rangle \langle (l+1)2^{j}+1\rangle } &
\cdots & P^{\langle (l+1)2^{j}+0\rangle \langle (l+1)2^{j}+(2^j-1)\rangle }\\
P^{\langle (l+1)2^{j}+1\rangle \langle (l+1)2^{j}+0\rangle }& P^{\langle (l+1)2^{j}+1\rangle \langle (l+1)2^{j}+1\rangle } &
\cdots & P^{\langle (l+1)2^{j}+1\rangle \langle (l+1)2^{j}+(2^j-1)\rangle }\\
\vdots & \vdots & \ddots &\vdots \\
P^{\langle (l+1)2^{j}+(2^j-1)\rangle \langle (l+1)2^{j}+0\rangle }& P^{\langle (l+1)2^{j}+(2^j-1)\rangle 
 \langle (l+1)2^{j}+1\rangle } & \cdots & P^{\langle (l+1)2^{j}+(2^j-1)\rangle \langle (l+1)2^{j}+2^j-1\rangle }
 \end{matrix}\right]~, \nonumber
\end{eqnarray}
and the inversion of the $X^l$ is,
\begin{equation}
\left(X^{l}_{[2^{j+1} \times 2^{j+1}]}\right)^{-1} =
\left[ \begin{matrix} \left(S_D^{\langle l\rangle }\right)^{-1}
& -\left(A^{l}\right)^{-1}B^l \left(S_A^{\langle l\rangle }\right)^{-1}\\
  -\left(D^{l}\right)^{-1}C^l \left(S_D^{\langle l\rangle }\right)^{-1}
  & \left(S_A^{\langle l\rangle }\right)^{-1}
\end{matrix} \right]~.
\label{eq:adinvprl}
\end{equation}
In the description of $X^l$, the Schur complements are computed as,
\begin{eqnarray}
S_D^{\langle l\rangle } & = &  A^l - B^l \left(D^l\right)^{-1} C^l~, \nonumber \\
S_A^{\langle l\rangle } & = &  D^l - C^l \left(A^l\right)^{-1} B^l~.
\end{eqnarray}
The list of intermediate calculations required during the process is
given in Table~\ref{tab:symbols}. We use suitable symbols to represent
the appropriate calculation for simplification in this article.
\begin{table}[ht]
\begin{center}
\begin{tabular}{|c|c|c|}
\hline
Symbol & Calculation & Description \\
\hline
$\Diamond$ & $\left(A^{l}\right)^{-1}$ and $\left(D^{l}\right)^{-1}$
& Inversion of diagonal blocks of $X^l$ \\ 
& $\left(S_A\right)^{-1}$ and $\left(S_D\right)^{-1}$ &
Inversion of Schur complements \\ \hline
$\rightarrow$ & $-\left(A^{l}\right)^{-1}B^l$ &
\multirow{2}{0.5\textwidth}{Multiplying the off diagonal blocks
with the inverted diagonal blocks.} \\ \cline{1-2}
$\leftarrow$ & $-\left(D^{l}\right)^{-1}C^l$ & \\ \cline{1-2}
\hline
$\sim$ & $S_A^{<l>}$ and $S_D^{<l>}$& Computing Schur complements \\
\hline
$\uparrow$ & $\left(-\left(A^{l}\right)^{-1}B^l\right)
  \left(S_A^{<l>}\right)^{-1}$ &
  \multirow{2}{0.5\textwidth}{Multiplying the off diagonal blocks
  with the inverted Schur complements.} \\ \cline{1-2}
$\downarrow$ & $\left(-\left(D^{l}\right)^{-1}C^l\right)
  \left(S_D^{<l>}\right)^{-1}$ & \\ \cline{1-2} \hline
\end{tabular}
\caption{Symbols to represent the calculations in the algorithm
described in Section \ref{ssec:inversexl}.}
\label{tab:symbols}
\end{center}
\end{table}
\subsection{Treatment of additional memory requirement and map
for computations in the algorithm}
\label{ssec:memory_parallel}
\begin{figure}
\begin{center}
\begin{tabular}{c c c c c c c c c c}
j & 0 & 0 & 0 & 0 & 0 & \dots & 0 & 0 & 0 \\
j & 1 & 0 & 0 & 0 & 0 & \dots & 0 & 0 & 0 \\
j & 2 & 0 & 0 & 0 & 0 & \dots & 0 & 0 & 0 \\
j & 2 & 1 & 0 & 0 & 0 & \dots & 0 & 0 & 0 \\
j & 3 & 0 & 0 & 0 & 0 & \dots & 0 & 0 & 0 \\
j & 3 & 1 & 0 & 0 & 0 & \dots & 0 & 0 & 0 \\
\vdots & \vdots & \vdots & \vdots & \vdots & \vdots & $\ddots$ & \vdots & \vdots & \vdots \\
j & j-1 & j-2 & j-3 & j-4 & j-5 & \dots & 3 & 1 & 0 \\
j & j-1 & j-2 & j-3 & j-4 & j-5 & \dots & 3 & 2 & 0 \\
j & j-1 & j-2 & j-3 & j-4 & j-5 & \dots & 3 & 2 & 1 \\
\end{tabular}
\caption{Array elements of {\sf loopid} in the steps to find
inverse of sub-block $X^l_{[2^{j+1}\times2^{j+1}]}$.}
\label{fig:loopid}
\end{center}
\end{figure} 
While finding inverses of sub-blocks $X^l_{[2^{j+1}\times2^{j+1}]}$,
as shown in Eq.~\ref{eq:adinvprl}, the computation of inverses of Schur
complements are to be performed iteratively. The inverse computations
for the further partitioned Schur complement blocks with lower orders,
i.e., sub-blocks with blockorder $[2^{x}\times2^{x}]$, in which $x\le j$,
cannot use the same memories of input matrix. This can be addressed by
creating additional memories before the beginning of full computation
process. By creating the additional memories as combination of set of
blocks, the additional advantage of efficient cache handling during
multiplication of sub-blocks can be achieved. For this purpose, set of
temporary matrices are created such that each set contains,

\begin{itemize}
\item $L$ - Left blocks (to store the computed values of $\leftarrow$), 
\item $R$ - Right blocks (to store the computed values of $\rightarrow$), 
\item $S$ - Square blocks (to store the computed values of $\sim$).
\end{itemize}
The number of provisional matrices is given by, $T=[({\rm
int})\log_2 (m_n)]-1$. Here the $L$ and $R$ blocks in the provisional matrices have
$k/2$ rows and the $S$ blocks have $k$ rows. The number of columns in the $T_l$'s
are $2^{l-1}$. Each block might have different order and they are fixed
according to the partitioning scheme of the input matrix so that the $T_{i}^{th}$
provisional matrix set blocks have {\sf blockorder} of $[2^i \times 2^i]$. The
computation is performed in distinctive steps in which each step is
identified by the variable `{\sf stepid}'. A map is created for each step which
is needed to be carried out distinctively, as a set of numbers which
we call as `{\sf loopid}'. Here it is to be noted that `{\sf loopid}' is an array
consists of $n$ elements, where $n$ is given by,
\begin{equation} \nonumber
n = \log_2~{\sf blocksize}+1~.
\end{equation}

Here, `{\sf blocksize}' represents the number of diagonal blocks. The
elements of `{\sf loopid}' are fixed based on `{\sf stepid}' and its
structure can be interpreted from Fig.~\ref{fig:loopid} for finding out
the inverses of sub-blocks $X^l_{[2^{j+1}\times2^{j+1}]}$ given that
inverses of $X^m_{[2^{j}\times2^{j}]},~m\in {1, 2, \dots (2l-1), 2l}$
is already known. The value of the element in the $i^{th}$ location for
a given {\sf stepid} follows a pattern and it cannot exceed the value
at $(i-1)^{th}$ location. For example, if there is a value $x$ in the
$i^{th}$ location in a particular step for $i<l$, then the value $x$
is retained in the upcoming steps, until the values at the locations
which are after $i$ form a decrement series by 1, i.e., the value at
the $(i+1)^{th}$ location is $x-1$, that at the $(i+2)^{th}$ location
is $x-2$,, etc., with the value at the $l^{th}$ location being 1 and
zero at the locations $i>l$. The above pattern can be recognized from
Fig.~\ref{fig:loopid}. This map is used for three mandatory purposes
which are,
\begin{enumerate}

\item Identification of what type of process is calculated among
$\leftarrow$, $\rightarrow$, $\sim$, $\uparrow$ and $\downarrow$ in
the given step. From the last non-zero element in the `{\sf loopid}'
from the left, the type of calculation to be carried out is decided.

\item To identify the sub-blocks $A_l$, $B_l$, $C_l$ and $D_l$ for the
step. From the location of last non-zero element in the `{\sf loopid}',
the presence of sub-blocks either of the input matrix or correct
provisional matrix set, which are required for the current calculations,
are identified.

\item To identify the set in which the calculated values are stored during
the calculation of $\leftarrow$, $\rightarrow$ and $\sim$. If the process
is meant to compute the Schur complement for a certain order, then the
last value of the non-zero element in the `{\sf loopid}' provides the
details of the provisional matrix set where the calculated values are
needed to be stored.

\end{enumerate}
Details of the procedure and the role of {\sf loopid} are given in the
next section.

\subsection{Complete procedure of large partitioned matrix inversion}
\label{ssec:complete}
\begin{table}[bth]
\centering
\begin{tabular}{|c|c c c c c c c c c|}
\hline
{\sf stepid} & \multicolumn{9}{c|}{{\sf loopid}}\\
\hline
1 & 1 & 0 & 0 & 0 & 0 & 0 & $\ldots$ & 0 & 0 \\
2 & 2 & 0 & 0 & 0 & 0 & 0 & $\ldots$ & 0 & 0 \\
3 & 2 & 1 & 0 & 0 & 0 & 0 & $\ldots$ & 0 & 0 \\
4 & 3 & 0 & 0 & 0 & 0 & 0 & $\ldots$ & 0 & 0 \\
5 & 3 & 1 & 0 & 0 & 0 & 0 & $\ldots$ & 0 & 0 \\
6 & 3 & 2 & 0 & 0 & 0 & 0 & $\ldots$ & 0 & 0 \\
7 & 3 & 2 & 1 & 0 & 0 & 0 & $\ldots$ & 0 & 0 \\
8 & 4 & 0 & 0 & 0 & 0 & 0 & $\ldots$ & 0 & 0 \\
9 & 4 & 1 & 0 & 0 & 0 & 0 & $\ldots$ & 0 & 0 \\
10 & 4 & 2 & 0 & 0 & 0 & 0 & $\ldots$ & 0 & 0 \\
11 & 4 & 2 & 1 & 0 & 0 & 0 & $\ldots$ & 0 & 0 \\
12 & 4 & 3 & 0 & 0 & 0 & 0 & $\ldots$ & 0 & 0 \\
\vdots & \multicolumn{9}{c|}{\vdots} \\
$N_{s} -2$ & $n$ & $n-1$ & $n-2$ & $n-3$ & $n-4$ & $n-5$ & $\dots$ & 1 & 0 \\
$N_{s} -1$ & $n$ & $n-1$ & $n-2$ & $n-3$ & $n-4$ & $n-5$ & $\dots$ & 2 & 0 \\
$N_{s}$ & $n$ & $n-1$ & $n-2$ & $n-3$ & $n-4$ & $n-5$ & $\dots$ & 2 & 1 \\
\hline
\end{tabular}
\caption{Mapping of the steps and their corresponding {\sf
loopid} in the algorithm. Here $N_{s}=(2*{\sf blocksize}-1)$ and $n$ is
the size of the {\sf loopid} array which is $n=\log_2~{\sf blocksize} +1$.}
\label{tab:map}
\end{table}
We divide the full calculation of inverse of partitioned matrix
into distinctive steps in which the processes can be performed in
parallel. Each step has to be executed and completed by
all the threads in full before moving to the next step which requires
an implicit thread barrier. For the above purpose we specifically
control the thread barriers with a key variable '{\sf stepid}'. For
a matrix partitioned into $k$ blocks, the number of {\sf stepid}'s are
$N_s=(2*{\sf blocksize}-1)$. We declare an array of {\sf loopid} with number of
elements as $n=\log_2~{\sf blocksize} +1$. The array elements are fixed
based on the {\sf stepid}. The calculations which are performed in each
step are decided based on the elements in the {\sf loopid} array. The
pattern in which the {\sf loopid} is fixed for different steps is provided
in Table~\ref{tab:map}.

The calculations involved to compute the inverses have already been
discussed and are listed in Table~\ref{tab:symbols}. The computations
which need to be performed, the location of input sub-blocks and the
storing place of calculation involved in each step are identified by the
`{\sf loopid}' as follows:

\begin{itemize}
\item If the `{\sf loopid}' contains first element as $1$ and remaining
elements are zero, it indicates the beginning of the inversion calculation
in which the inversion of diagonal blocks in input partitioned matrix
of blockorder $[1 \times 1]$, are performed and stored at the inverted
matrix.

\item If the first element in the `{\sf loopid}' is $j$ where ($j>1$)
and all other elements are zero, then the processes which need to be
performed are the calculation of $\leftarrow$, $\rightarrow$ and $\sim$s
of the input partitioned matrix. Hence, the sub-blocks required in the
calculations--- $A_l$, $B_l$, $C_l$ and $D_l$ ---are identified in the
input matrix and each sub-block has the blockorder $\left[ 2^{j-1} \times
2^{j-1} \right]$. The calculated values are stored in the $(j-1)^{th}$
provisional matrix set which is identified from the value of $j$.

\item If the first non-zero element of the `{\sf loopid}' is $1$
and occurs as the $x^{th}$ element, $x \ne 1$, then the inversion
of diagonal blocks present in the square blocks ($S$) of provisional
matrix set $T_{cid}$, need to be performed. Here the provisional matrix
set, $T_{cid}$, is identified from the value of element, say $p_1$
at $(x-1)^{th}$ location in the `{\sf loopid}' and then the $id$ is
$p_1-1$. The elements of `{\sf loopid}' for this process are arranged as,
\begin{center}
\begin{tabular}{c | c c c c c c c c c} 
Location in the `{\sf loopid}' & 1 & \ldots & $x-3$ & $x-2$ &
		$x-1$ & x & $x+1$ & $x+2$ & \ldots\\ 
Elements of `{\sf loopid}' & $j$ & \ldots & $p_3$ & $p_2$ &
$p_1$ & 1 & 0 & 0 & \ldots \\
\end{tabular}
\end{center}
This step calculates and stores the inverse of diagonal blocks
({\sf blockorder} $[2^{1} \times 2^{1}]$) of Schur complements in
the inverted matrix. If $p_1$ is $2$, then the inverse for blockorder
$[2^{2} \times 2^{2}]$, can be calculated as $\uparrow$ and $\downarrow$
using the calculated inverses for the provisional matrix set $T_{cid}$
by treating the non-diagonal blocks as $B_l$s and $C_l$s in the square
blocks of the set $S$. It results in inverse of partitioned blocks of
blockorder $[2^{2} \times 2^{2}]$. The process of calculating $\uparrow$
and $\downarrow$ is iteratively carried forward if $p_{i+1} = p_{i}$,
and the result will yield inverted blocks of block order $[2^{i+1}
\times 2^{i+1}]$ in the the inverted matrix.

\item If the first non-zero element of the `{\sf loopid}' is $z$, $z>1$
at the location $x$, then the process is to compute the $\leftarrow$,
$\rightarrow$ and $\sim$ using the square blocks $S$ present in the
provisional matrix set $T_{cid}$ and the calculated values are stored
at the respective blocks $L$, $R$ and $S$ of provisional matrix set
$T_{sid}$. Here $T_{cid}$, is identified from the value of element, say
$p_1$ at $(x-1)^{th}$ location in the `{\sf loopid}' and then the $cid$
is $p_1-1$ and the $T_{sid}$ where the calculated values are stored is
identified from the value $z$ which is at the location $x$, as $sid=z-1$.
\end{itemize}

The calculations are such that they have to be completed by all the
threads in each step corresponding to {\sf stepid} before executing
the next, since subsequent calculations would depend on the previous
execution. The last non-zero value and its location in the {\sf loopid}
decides what calculation has to be performed.  The pattern for the values {\sf loopid} over different {\sf stepid} can be recognized from Table~\ref{tab:map}.

It is to be noted that the complete {\sf loopid} elements for any {\sf stepid} can be fixed only with the value of the {\sf stepid}. It means that the {\sf stepid} serve as a control variable for applying thread barriers which are required to control the calculations in which the calculations in each step has to be carried out after completely executing the calculations in previous steps among different threads during the parallel processing. In this context {\sf stepid} is a crucial control variable in this algorithm.
\begin{center}
\begin{table}[]
\begin{tabular}{|p{0.07\textwidth}|p{0.095\textwidth}|p{0.095\textwidth}|
p{0.12\textwidth}|p{0.12\textwidth}|p{0.12\textwidth}|p{0.2\textwidth}|} \hline
{\sf stepid} & Last non-zero element in {\sf loopid} & Location of
last non-zero element in {\sf loopid} & Calculation & Calculated using
the blocks in & Place of storing or updating & {\sf loopid} \\ \hline
1&1&1& $\Diamond$ & $M$ & $M_{inv}$ & $\underline{1}~0~0~0~0~...~0~0~0$ \\
\hline
2&2&1& $\rightarrow ~\leftarrow ~\sim $ & $M$ & $T_{1}$ & $\underline{2}~0~0~0~0~...~0~0~0$ \\
\hline
3&1&2& $\Diamond$ & $T_1$ &  \multirow{2}{3cm}{$M_{inv}$} & {$2~\underline{1}~0~0~0~...~0~0~0$} \\
&&& $\uparrow~\downarrow$ &$T_1$&& {$\underline{2~1}~0~0~0~...~0~0~0$}\\
\hline
4&3&1& $\rightarrow ~\leftarrow ~\sim $ & $M$ & $T_{2}$ & $\underline{3}~0~0~0~0~...~0~0~0$ \\
\hline
5&1&2& $\Diamond$ & $M$ & $M_{inv}$ & $3~\underline{1}~0~0~0~...~0~0~0$ \\
\hline
6&2&1& $\rightarrow ~\leftarrow ~\sim $ & $M$ & $T_{1}$ & $3~\underline{2}~0~0~0~...~0~0~0$ \\
\hline
7&1&2& $\Diamond$ & $T_1$ &  \multirow{3}{3cm}{$M_{inv}$} & {$3~2~\underline{1}~0~0~...~0~0~0$} \\
&&& $\uparrow~\downarrow$ &$T_1$&& {$3~\underline{2~1}~0~0~...~0~0~0$}\\
&&& $\uparrow~\downarrow$ &$T_2$&& {$\underline{3~2}~1~0~0~...~0~0~0$}\\
\hline
8&4&1& $\rightarrow ~\leftarrow ~\sim $ & $M$ & $T_{3}$ & $\underline{4}~0~0~0~0~...~0~0~0$ \\
\hline
9&1&2& $\Diamond$ & $M$ & $M_{inv}$ & $4~\underline{1}~0~0~0~...~0~0~0$ \\
\hline
10&2&2& $\rightarrow ~\leftarrow ~\sim $ & $M$ & $T_{1}$ & $4~\underline{2}~0~0~0~...~0~0~0$ \\
\hline
11&1&3& $\Diamond$ & $T_1$ &  \multirow{2}{3cm}{$M_{inv}$} & {$4~2~\underline{1}~0~0~...~0~0~0$} \\
&&& $\uparrow~\downarrow$ &$T_1$&& {$4~\underline{2~1}~0~0~...~0~0~0$}\\
\hline
\vdots & \vdots & \vdots &\vdots &\vdots &\vdots &\vdots \\
\hline
$N_{s}-2$&1&$n-1$& $\Diamond$ & $M$ & $M_{inv}$ & {$n~n-1~..........~3~\underline{1}~0$} \\
\hline
$N_{s}-1$&2&$n-1$& $\rightarrow ~\leftarrow ~\sim $ & $M$ & $T_{1}$ & {$n~n-1~..........~3~\underline{2}~0$} \\
\hline
$N_{s}$&1&$n$& $\Diamond$ & $T_1$ &  \multirow{5}{3cm}{$M_{inv}$} & {$n~n-1~..........~3~2~\underline{1}$} \\
&&& $\uparrow~\downarrow$ &$T_1$&& {$n~n-1~..........~3~\underline{2~1}$}\\
&&& $\uparrow~\downarrow$ &$T_2$&& {$n~n-1~..........~\underline{3~2}~1$}\\
&&& \vdots &\vdots && \vdots\\
&&& $\uparrow~\downarrow$ &$T_{n-1}$&& {$\underline{n~n-1}~........~3~2~1$}\\
\hline
\end{tabular}
\caption{Complete procedure to invert the matrix which is partitioned in
blocks form of block order $[N_k \times N_k]$. The calculations performed in each step are identified from the underlined elements shown in the last column and the location of those elements at `{\sf loopid}'.}
\label{tab:complete}
\end{table}
\end{center}

In Table~\ref{tab:complete}, the complete procedure for inverting
the matrix is given. In each step, the processes can be performed
in parallel among multiple threads. At any instance of the calculation, the
procedure is such that the calculation never accesses the same memory
blocks between the threads for storing and also the input blocks in
the calculation too. The barriers are controlled using the {\sf stepid}
which will ensure the recursion as bottom to top approach. The other
requirement which is freedom for the threads to select a part of a
large operation is also being satisfied here. And once the parallel
process is started, the algorithm does not require closing and reopening
the parallel process.

\subsection{Using the Fox method of matrix multiplication in the
$\leftarrow$ and $\rightarrow$ calculations}
\label{ssec:fox}
Large matrix multiplications can be performed using Fox's algorithm
for parallel processing~\cite{divideblockmul, fox, scaling}. Usually,
the large matrices are blocked in such a way that all the blocks are
of the same order for implementing this procedure. In this method,
every iteration is performed with efficient cache handling as divide
and conquer way \cite{cachemul}. In our procedure, the input matrix is
already partitioned into blocks such that the diagonal blocks are square
matrices. And for the calculations of $\leftarrow$ and $\rightarrow$, it
can be seen from the description of $A^l, B^l, C^l$ and $D^l$ that the
row-orders of partitioned blocks in the $A^l$ blocks and $D^l$ blocks
match the column-orders of blocks $B^l$ and $C^l$ respectively. Also,
the orders of partitioned blocks in $L$ and $R$ in the $T$ provisional
matrix set, in which the computed $\leftarrow$ and $\rightarrow$ will be
stored, ensures this multiplication can be performed blockwise. Hence,
Fox's method for multiplications among the sub-blocks can be implemented
even though all the sub-blocks are not of the same order.

Since our partitioning scheme can be varied keeping only the number of
diagonal blocks as powers of 2, we have the advantage of optimizing
for cache handling in the calculation. As an added advantage these
multiplications can also be performed in parallel without
altering the partitioned form for every blockorder in $X^l$.
\subsection{An efficient method for performing $\uparrow$ and $\downarrow$
calculations}
\label{ssec:arrowops}
In Table~\ref{tab:complete}, the calculation of $\uparrow$ and
$\downarrow$ are performed in the steps $3, 7, 11, ... , N_s$ after
completing the inverse of diagonal blocks. At a given particular step,
every $\uparrow$ and $\downarrow$ calculation provides the inverses
of blocks with blockorder as one additional power of 2 in inverted
matrix each time. The structure of {\sf loopid} which perform $x$
times the $\uparrow$ and $\downarrow$ calculations in a step is given by
\begin{center}
$\cdots ~y~x~(x-1)~(x-2)~\cdots ~3~2~1~ \cdots ~$.
\end{center}

After computing the inverses of the diagonal blocks, $\uparrow$ and
$\downarrow$ calculations have to performed as follows:

\begin{itemize}
\item $\cdots y~x~(x-1)~(x-2)~\cdots ~3~\underline{2~1}~0~0~ \cdots~$:
This step describes the right multiplication of inverted blocks in the
$M_{inv}$ of order $\left[ 2^1 \times 2^1 \right]$ which result from
the $\Diamond$ calculation, with the $R$ blocks in $T_{1}$. Since the
calculated values are stored in $M_{inv}$, this step provides the inverse
of blocks of order $\left[ 2^2 \times 2^2 \right]$.
\item $\cdots y~x~(x-1)~(x-2)~\cdots ~\underline{3~2}~1~0~0~ \cdots~$:
This step describes the right multiplication of inverted blocks in
the $M_{inv}$ of order $\left[ 2^2 \times 2^2 \right]$ with the $R$
blocks in $T_{2}$ and the calculated values are used to further update
$M_{inv}$ with the calculated inverses of blocks of order $\left[ 2^3
\times 2^3 \right]$.

{\hspace*{7cm} \textbf{$\vdots$}}
\item $\cdots y~\underline{x~(x-1)}~(x-2)~\cdots~3~2~1~0~0~\cdots~$:
This step describes the right multiplication of inverted blocks in the
$M_{inv}$ of order $\left[ 2^{x-1} \times 2^{x-2} \right]$ with the $R$
blocks in $T_{x-1}$ and the calculated values are used to further update
$M_{inv}$ with the calculated inverses of blocks of order $\left[ 2^x
\times 2^x \right]$.
\end{itemize}
It can be observed that starting with inverses of diagonal blocks of
order $\left[ 2^{1} \times 2^{1} \right]$ in the $M_{inv}$, the inverses
of blocks of order $\left[ 2^x \times 2^x \right]$ are obtained at the
end of the calculation in this step. Due to the fact that the diagonal
blocks are square matrices in $M_{inv}$, order of the blocks in every
partitioned $X^l$, ($l\ge 1$) row-wise, will match with the order of
blocks in $L$ and $R$ for given $T^l$. In right multiplication in the
$\uparrow$ and $\downarrow$, it can be seen that partitioned blocks are
iteratively being multiplied with the $L$ and $R$ blocks of different
$T^i$ provisional matrix set. The point here is all different $\uparrow$
and $\downarrow$ can be clubbed together and performed iteratively
to obtain the inverses of blocks of order $\left[ 2^x \times 2^x
\right]$. For a sub-block $X^l$, the number of repetitive calculations
which has to be performed on the sub-sub-blocks of $M_{inv}$ to complete
all the $\uparrow$ and $\downarrow$ calculations in a step can be mapped
as in Fig.~\ref{fig:updownmap}.
\begin{figure}[htp]
\begin{equation*}
\left( \begin{matrix}
1 & 2 & 2 & 3 & 2 & 3 & 3 & 4 & \cdots ~\\
2 & 1 & 3 & 2 & 3 & 2 & 4 & 3 & \cdots ~\\
2 & 3 & 1 & 2 & 3 & 4 & 2 & 3 & \cdots ~\\
3 & 2 & 2 & 1 & 4 & 3 & 3 & 2 & \cdots ~\\
2 & 3 & 3 & 4 & 1 & 2 & 2 & 3 & \cdots ~\\
3 & 2 & 4 & 3 & 2 & 1 & 3 & 2 & \cdots ~\\
3 & 4 & 2 & 3 & 2 & 3 & 1 & 2 & \cdots ~\\
4 & 3 & 3 & 2 & 3 & 2 & 2 & 1 & \cdots ~\\
\vdots & \vdots & \vdots & \vdots & \vdots & \vdots & \vdots & \vdots & \ddots \\
\end{matrix} \right)
\end{equation*}
\caption{Map of performing $\downarrow$ and $\uparrow$ calculations in a
block $X^l$ required by the steps which have last non-zero element
in the {\sf loopid} as 1, in the Table \ref{tab:complete}.}
\label{fig:updownmap}
\end{figure}
This method starts with the computation of inversion of {\sf blockorder}
$\left[ 2^0 \times 2^0 \right]$ on $M_{inv}$. It is represented as $1$
(in the diagonal blocks) in Fig.~\ref{fig:updownmap}. In the second
iteration, the $\left[ 2^0 \times 2^0 \right]$ blocks are right multiplied
with the remaining blocks in the same column appropriately with the blocks
in the provisional matrices set. It results in the inversion of blockorder
$\left[ 2^1 \times 2^1 \right]$ on $M_{inv}$ along with partially computed
values in the remaining columns of $X^l$. For the next iteration-3,
the completely inverted sub-blocks other than the diagonal in $\left[
2^1 \times 2^1 \right]$ scheme, are loaded into the memory and right
multiplied appropriately with the blocks in provisional matrix set, and
added with the values already present in all the remaining blocks which
are not yet a part of the full inversion by reduction operation. This
procedure is further repeated for obtaining inverses of sub-block of next
power of 2. Each time the reduction is performed over the other blocks
which only provides partially computed values in the blocks. The number
of iterations required to completely obtain the inverse in each sub-block
is given as map in Fig.~\ref{fig:updownmap}.

The procedure shows that we can avoid reloading of different blocks
of $M_{inv}$ for different sets of calculations of $\uparrow$ and
$\downarrow$. In this map parallelization is performed column-wise
since the calculations in every column only depend on the previously
calculated values in the same column. For every iteration, calculations
done in the previous iteration of the calculation are required to compute
forward. As given above this combined calculation reduces the repetitive
loading of same blocks of $M_{inv}$ into the cache.

\subsection{Discussion about the performance and significance of the algorithm}
\label{ssec:performance}
\begin{center}
\begin{figure}[htb]
\centering
\includegraphics[scale=0.49, angle=270]{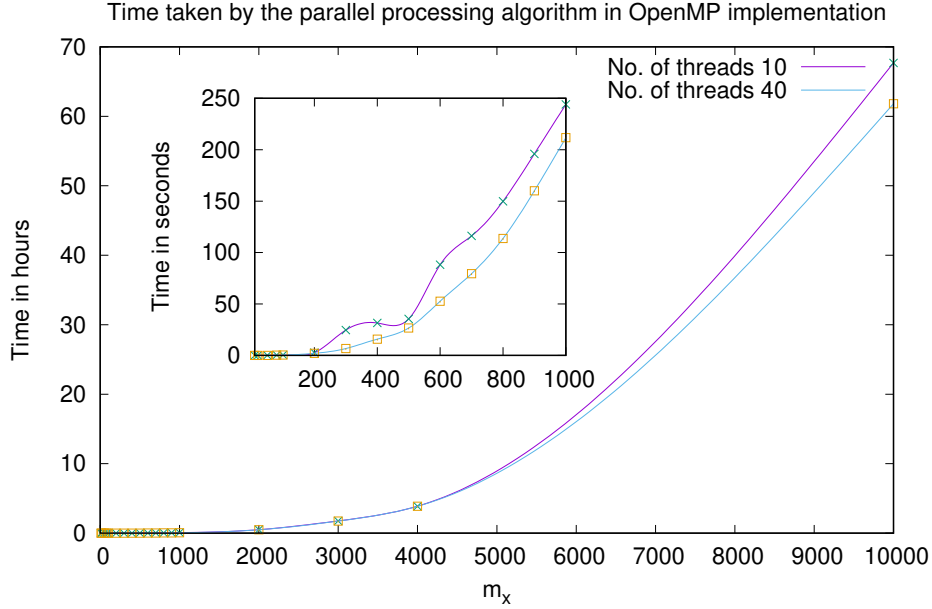}
\caption{Time taken for inverting matrices up to $m_x=10^4$ using the OpenMP implementation of developed algorithm in Nandadevi cluster for the number of threads 10 and 40.  The inside figure shows the computation time for matrices of order up to 1000.}
\label{fig:prlltime}
\end{figure}
\end{center}
\begin{center}
\begin{figure}[htb]
\centering
\includegraphics[scale=0.4568, angle=270]{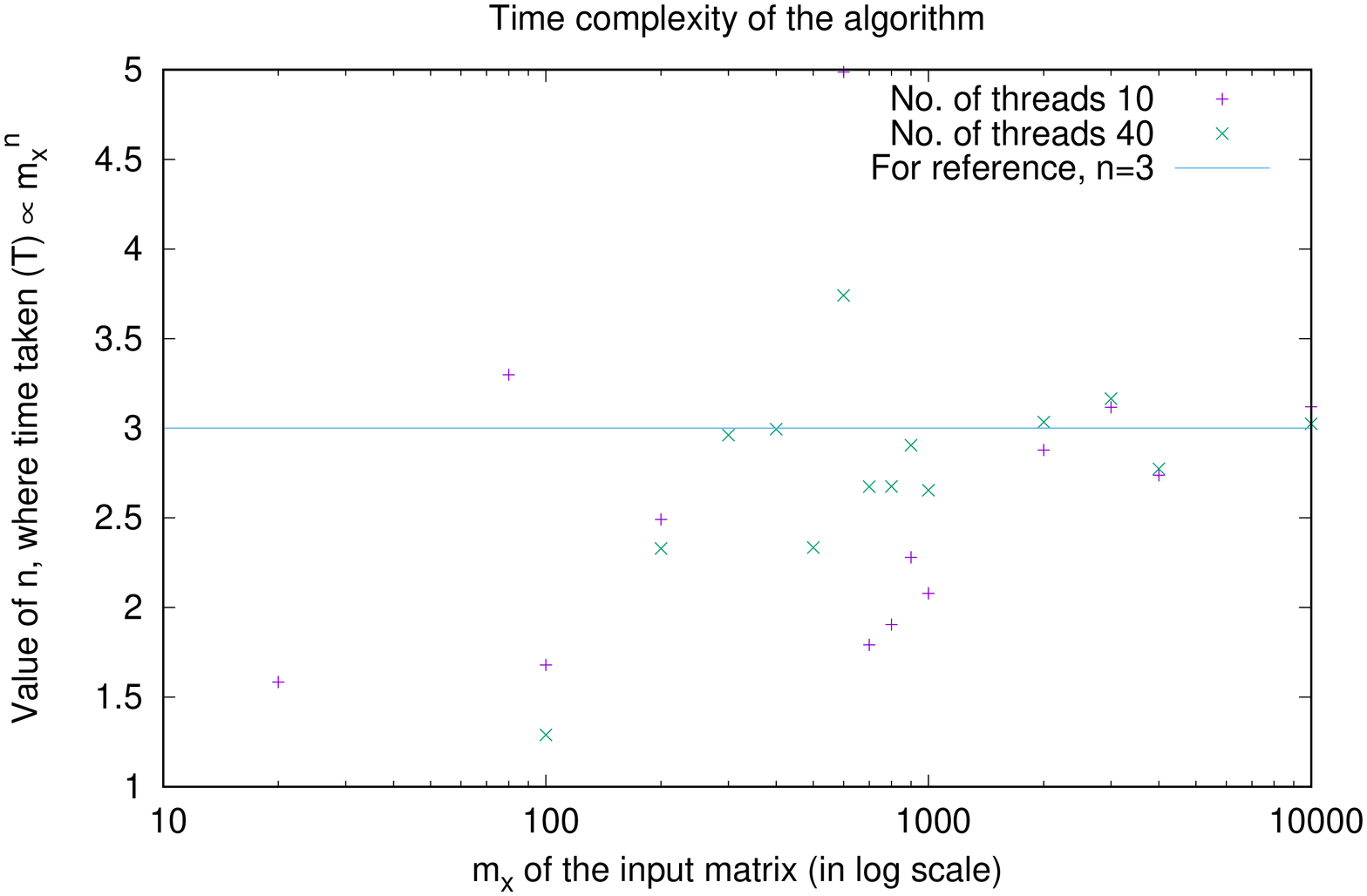}
\caption{Time complexity of algorithm from the OpenMP implementation in Nandadevi cluster.  With the assumption for time taken for computing matrix inversion ($T$) as, $T \propto m_x^n$, the value of $n$ in the time complexity $\mathcal{O}\left(m_x^n\right)$, is calculated from the derivative of $log~T$ with respect to $log~m_x$ at different values of $m_x$.  The figure shows the value of n in time complexity fluctuates around 3 for the matrices with $m_x$ up to $10^4$.}
\label{fig:prllslope}
\end{figure}
\end{center}
The developed algorithm was implemented using OpenMP application
programming interface (API) and written in C language. The code was tested
in Nandadevi cluster for parallel processing performance at the Institute
of Mathematical Sciences with the random generated matrices up to the order
$m_x \approx 10^4$. The node environment in the cluster used in the analysis comprises of two 10C Intel Xeon E5-2687W (3.1 GHz) processors and 64GB DDR3 1866MHz RAM. The OpenMP version of parallel processing for matrix inversion is not expected to improve significantly in the timing for inversion calculation than serial processing, because
of the fact that matrix inversions calculations are highly dependent on the shared memory handling. In Fig. \ref{fig:prlltime}, the timing performance of algorithm for matrices with $m_x$ up to $10^4$ is shown for computing the inversion with 10 and 40 threads.  In Fig. \ref{fig:prllslope}, the calculated value of $n$, in the time complexity, for different values of $m_x$ is shown.  The value of $n$ fluctuates around $3$ most of the time.  The calculated value of $n$ for 10 and 40 threads are $2.35 \pm 0.13$ and $3.01 \pm 0.006$ for the $m_x$ range of $10$ to $10^4$.  The observation is that even within the OpenMP implementation, we infer the timing complexity is in the order of $m_x^3$ approximately.
The code and its performance which are presented here is to establish the capability
of this developed method to invert general matrices by partitioning into
large number of blocks. The same procedure can also be in principle applied to do the computation with Graphics processing unit (GPU), where the improvement
in timing for parallel processing could be seen, as GPUs are the better place for implementing parallel methods when the calculations are highly dependent on the ways in which memories are handled.
It is to be noted that the choice of orders in the diagonal blocks during the partitioning is completely arbitrary.  Here, the code is presented for a simple scenario of sub diagonal blocks order are chose from the set ${2, 3, 4}$.  But, in realistic situation, the choice could be as large as the computer can handle.  And the method of inverting the partitioned sub blocks are also arbitrary.  The other catch is that the multiplications which are required during the inversions are performed in same partitioned blocks structure format, by updating the blocks with new calculated values when and where required.  We have observed that even for writing the input matrix and inverted matrix with $m_x=7000$ into a file take nearly $1~$GB of storage space when the elements are treated as doubles.  But, we can have scenario where the partitioned blocks are stored in separate files.   Since the required amount of memory including for the storage of temporary memory blocks is defined, handling the situation of large partitioned blocks and the blocks of provisional matrices among the separate files is possible.  We also considered the scenario where the calculations are carried out with many complete stops of the machine, since the calculations do not have to perform in one go fully, but can be performed in a predefined sequence with appropriate breaks using the stepid provided in the parallel algorithm.  It means that we may overcome both memory and time limitations for performing inverses of extremely large matrices.

Gauss Jordan Elimination method and LU decomposition method are the widely used parallel computing methods for calculating the matrix inversion and there are techniques to compute matrix inversion for several special types of matrices \cite{ultra, prll1, prll2, prll3}.  But, the algorithm which is presented in this article is meant to perform the computation for different cases of block partitioned scheme in which the order of the partitioned blocks can be arbitrarily chosen.  Also, this algorithm sequentially combine the inverted blocks to the next blockorder in the powers of two.  Hence, it is significant and the article established the possibility of performing matrix inversion using parallel processing with blockwise inversion technique itself.  For a greater improvement in terms of computation time, the algorithm has to be implemented in GPUs.  It can be seen that multi-dimensional parallellization is possible for calculating $\uparrow$ and $\downarrow$, in the context of the presented method itself.  Hence, for future study, we intend to utilize the APIs such as CUDA which allows higher dimensional parallelization for the improvement in the presented algorithm.

\section{Conclusion}
\label{sec:concl}
The block structure are quite often encountered in different physics problems and one such incident of calculating systematic uncertainties for tau neutrino events in our analysis~\cite{taupaper} is addressed by calculating matrix inversion using blockwise inversion technique.  From further studies with blockwise inversion for different ways of memory handling, we addressed various scenarios of inverting the block matrix which are applicable in different fields of research.  This work provides the developed programs in C language for practical applications as `invertor' project.  The performance of developed programs with a common present day computer is presented in this article. A comparison study of different ways of block matrix inversion was made due to the difference in memory management and computational complexity. Three different procedures for blockwise inversion of matrix are discussed by assuming the diagonal blocks and
their Schur complements are invertible in the block partitioned
scheme. Flowchart diagrams adopting these procedures are provided
in Figures ~\ref{fig:inversionfull},~\ref{fig:inplaceinversion}
and~\ref{fig:invertbyaandd}.  Each procedure has different advantage and the optimum selection of method can be chosen according to the requirement. These procedures and the programs to invert large matrices which are presented here, are not meant to establish any superiority of their techniques over existing well performing programs in different languages. Rather, by means of the studies with the developed procedures, we establish the possibility of creating an algorithm for dealing with ultra large matrices to overcome the different limitations due to memory management in computations. The intermediate steps in the developed algorithm performs the computation by treating the large partitioned blocks as if the blocks are just numbers and blocks structure is preserved throughout the calculation. The choice of blockorder for the blocks is completely arbitrary and the memory management for the calculations are precisely defined.  The restriction is at the level of
number of partitioned blocks in diagonal which need to be in the power
of 2. This algorithm is developed for performing block matrix inversion
in clusters using parallel processing which is an added advantage.  

\paragraph{Acknowledgments}: The author is thankful to D. Indumathi,
M.V.N. Murthy, Sitabhra Sinha and Suprabh Prakash for discussions and would like to
acknowledge the High Performance Computing (HPC) facility---the Nandadevi
cluster at The Institute of Mathematical Sciences.

\end{document}